\documentstyle[12pt, amsfonts]{amsart}
\begin{document}
\def\R{{\mathbb R}}
\def\Z{{\mathbb Z}}
\def\C{{\mathbb C}}
\newcommand{\trace}{\rm trace}
\newcommand{\Ex}{{\mathbb{E}}}
\newcommand{\Prob}{{\mathbb{P}}}
\newcommand{\E}{{\cal E}}
\newcommand{\F}{{\cal F}}
\newtheorem{df}{Definition}
\newtheorem{theorem}{Theorem}
\newtheorem{lemma}{Lemma}
\newtheorem{pr}{Proposition}
\newtheorem{co}{Corollary}
\newtheorem{problem}{Problem}
\def\n{\nu}
\def\sign{\mbox{ sign }}
\def\a{\alpha}
\def\N{{\mathbb N}}
\def\A{{\cal A}}
\def\L{{\cal L}}
\def\X{{\cal X}}
\def\F{{\cal F}}
\def\c{\bar{c}}
\def\v{\nu}
\def\d{\delta}
\def\diam{\mbox{\rm dim}}
\def\vol{\mbox{\rm Vol}}
\def\b{\beta}
\def\t{\theta}
\def\l{\lambda}
\def\e{\varepsilon}
\def\colon{{:}\;}
\def\pf{\noindent {\bf Proof :  \  }}
\def\endpf{ \begin{flushright}
$ \Box $ \\
\end{flushright}}

\title[Isomorphic Busemann-Petty problem]
{Isomorphic Busemann-Petty problem for sections of proportional dimensions}

\author{Alexander Koldobsky}

\address{Department of Mathematics\\ 
University of Missouri\\
Columbia, MO 65211}

\email{koldobskiya@@missouri.edu}

\begin{abstract}  The main result of this note is a solution to the isomorphic Busemann-Petty problem 
for sections of proportional
dimensions, as follows. Suppose that $0<\lambda<1,$ $k>\lambda n,$ and $K,L$ are
origin-symmetric convex bodies in $\R^n$ satisfying the inequalities
$$|K\cap H|\le |L\cap H|,\qquad \forall H\in Gr_{n-k},$$
where $Gr_{n-k}$ is the Grassmanian of $(n-k)$-dimensional subspaces of $\R^n,$
and $|K|$ stands for volume of proper dimension.
Then
$$|K|^{\frac{n-k}n}\le C^k \left(\sqrt{\frac{(1-\log \lambda)^3}{\lambda}}\right)^k|L|^{\frac{n-k}n},$$
where $C$ is an absolute constant.

\end{abstract}  
\maketitle

\section{Introduction} The {Busemann-Petty problem}, raised in 1956 in \cite{BP}, 
asks the following question.
Suppose that $K,L$ are origin-symmetric convex bodies in $\R^n$ so that the $(n-1)$-dimensional volume of
every central hyperplane section of $K$ is smaller than the same for $L,$ i.e.
\begin{equation}\label{bp-condition}
|K\cap \xi^\bot|\le |L\cap \xi^\bot|,\qquad \forall \xi\in S^{n-1}.
\end{equation}
Does it follow that the $n$-dimensional volume of $K$ is smaller than that of $L,$ i.e.
$$|K|\le |L|\ ?$$
Here $\xi^\bot =\{x\in \R^n:(x,\xi)=0\}$ is the central hyperplane perpendicular to $\xi,$ 
and $|K|$ stands for volume of proper dimension.
The problem was solved in the end of the 1990's as the result 
of a sequence of papers \cite{LR}, \cite{Ba2}, \cite{Gi}, \cite{Bo4}, 
\cite{L}, \cite{Pa}, \cite{G1}, \cite{G2}, \cite{Z1}, \cite{Z2}, \cite{K1}, \cite{K2}, \cite{Z3},
\cite{GKS} ; see \cite[p. 3]{K4} or \cite[p. 343]{G3} for details.
The answer is affirmative if $n\le 4$, and it is negative if $n\ge 5.$

The {lower dimensional Busemann-Petty problem} asks the same question for
sections of lower dimensions. Suppose $K,L$ are origin-symmetric convex bodies in $\R^n,$
and $1\le k \le n-1.$ Let $Gr_{n-k}$ be the Grassmanian of $(n-k)$-dimensional subspaces
of $\R^n,$ and suppose that
\begin{equation}\label{ldbp-condition}
|K\cap H|\le |L\cap H|,\qquad \forall H\in Gr_{n-k}.
\end{equation}
Does it follow that $|K|\le |L|?$ It was proved in \cite{BZ} (see also \cite{K3}, \cite[p.112]{K4}, \cite{RZ} and 
\cite{M} for different proofs) that the answer is negative if the dimension of sections $n-k>3.$
The problem is still open for two- and three-dimensional sections ($n-k=2,3,\ n\ge 5).$

Since the answer to the Busemann-Petty problem is negative in most dimensions, it makes sense to ask the
{isomorphic Busemann-Petty problem}, namely, does there exist an absolute 
constant $C$ such that inequalities (\ref{bp-condition}) imply $$|K|\le C\ |L|\ ?$$
If the answer to the isomorphic Busemann-Petty problem was affirmative, then by iteration there would exist an absolute 
constant $C$ such that for every $1\le k \le n-1$ 
\begin{equation}\label{ildbp}
|K|^{\frac{n-k}n}\le C^k |L|^{\frac{n-k}n}.
\end{equation}
However, the isomorphic Busemann-Petty problem is still open and equivalent to the slicing problem \cite{Bo1, Bo2, Ba1, MP}, 
another major open problem in convex geometry. The slicing problem asks
whether there exists an absolute constant $C$ so that for any origin-symmetric convex body $K$ in $\R^n$
of volume 1 there is a hyperplane section of $K$ whose $(n-1)$-dimensional volume is greater than $1/C.$
In other words, does there exist an absolute constant $C$ so that for any $n\in \N$ and any
origin-symmetric convex body $K$ in $\R^n$
\begin{equation} \label{hyper}
|K|^{\frac {n-1}n} \le C \max_{\xi \in S^{n-1}} |K\cap \xi^\bot|.
\end{equation}
The best current result $C\le O(n^{1/4})$ is due to Klartag \cite{Kl}, who
removed the  logarithmic term from an earlier estimate of Bourgain \cite{Bo3}.
We refer the reader to [BGVV] for the history and partial results.

Iterating (\ref{hyper}) one gets the lower dimensional slicing problem asking whether
the inequality
\begin{equation} \label{lowdimhyper}
|K|^{\frac {n-k}n} \le C^k \max_{H\in Gr_{n-k}} |K\cap H|
\end{equation}
holds with an absolute constant $C,$ where $1\le k \le n-1.$ 
Inequality (\ref{lowdimhyper}) was recently proved in \cite{K5} in the case
where $k\ge \lambda n,\ 0<\lambda<1,$ with the constant $C=C(\lambda)$ dependent only on $\lambda.$

\begin{pr}\label{ldsp} (\cite[Corollary 3]{K5}) There exists an absolute constant $C$ such that for every $n\in \N,$
every $0<\lambda<1,$ every $k>\lambda n,$ and every origin-symmetric convex body $K$ in $\R^n$
$$ |K|^{\frac{n-k}n}\ \le\  C^k 
\left(\sqrt{\frac{(1-\log \lambda)^3}{\lambda}}\right)^k \max_{H\in Gr_{n-k}} |K\cap H|.$$
\end{pr}

In this note we prove an isomorphic Busemann-Petty problem for sections of proportional
dimensions.

\begin{theorem} \label{ildbp-th} Suppose that $0<\lambda<1,$ $k>\lambda n,$ and $K,L$ are
origin-symmetric convex bodies in $\R^n$ satisfying the inequalities
$$|K\cap H|\le |L\cap H|,\qquad \forall H\in Gr_{n-k}.$$
Then
$$|K|^{\frac{n-k}n}\le C^k \left(\sqrt{\frac{(1-\log \lambda)^3}{\lambda}}\right)^k|L|^{\frac{n-k}n},$$
where $C$ is an absolute constant.
\end{theorem}
It is easy to see that Theorem \ref{ildbp-th} implies Proposition \ref{ldsp}; see Remark 2. It is not known whether
Theorem \ref{ildbp-th} can be deduced from Proposition \ref{ldsp}, we provide an independent proof
here. 

Proposition \ref{ldsp} was proved in \cite{K5} for arbitrary measures 
in place of volume. The arguments of this paper do not allow for an extension of Theorem \ref{ildbp-th}
to arbitrary measures, and, therefore, the possibility of such an extension remains open.
Note that  a version of the isomorphic Busemann-Petty problem for arbitrary measures was established in \cite{KZ},
but with the constant $C=\sqrt{n}$ depending on the dimension.

\section{Proof of Theorem \ref{ildbp-th}} \label{slicing}

We need several definitions and facts.
A closed bounded set $K$ in $\R^n$ is called a {\it star body}  if 
every straight line passing through the origin crosses the boundary of $K$ 
at exactly two points different from the origin, the origin is an interior point of $K,$
and the {\it Minkowski functional} 
of $K$ defined by 
$$\|x\|_K = \min\{a\ge 0:\ x\in aK\}$$
is a continuous function on $\R^n.$ 

We use the polar formula for volume of a star body
\begin{equation}\label{polar}
|K|=\frac 1n \int_{S^{n-1}} \|\theta\|_K^{-n} d\theta.
\end{equation}

The solution ot the original Busemann-Petty problem was based on a connection with intersection bodies 
found by Lutwak \cite{L}. In this paper we use a more general class of generalized  intersection bodies
introduced by Zhang \cite{Z4} 
in connection with the lower dimensional Busemann-Petty problem.

For $1\le k \le n-1,$  the {\it $(n-k)$-dimensional spherical Radon transform} 
$R_{n-k}:C(S^{n-1})\to C(Gr_{n-k})$  
is a linear operator defined by
$$R_{n-k}g (H)=\int_{S^{n-1}\cap H} g(x)\ dx,\quad \forall  H\in Gr_{n-k}$$
for every function $g\in C(S^{n-1}).$ By the polar formula for volume, for every
$H\in Gr_{n-k},$ we have
\begin{equation}\label{polar-ld}
|K\cap H| =\frac 1{n-k}\int_{S^{n-1}\cap H} \|x\|_K^{-n+k}dx = \frac 1{n-k} R_{n-k}(\|\cdot\|_K^{-n+k})(H).
\end{equation}

We say that an origin symmetric star body $K$ in $\R^n$ is a {\it generalized $k$-intersection body}, 
and write $K\in {\cal{BP}}_k^n,$  if there exists a finite Borel non-negative measure $\mu$
on $Gr_{n-k}$ so that for every $g\in C(S^{n-1})$
\begin{equation}\label{genint}
\int_{S^{n-1}} \|x\|_K^{-k} g(x)\ dx=\int_{Gr_{n-k}} R_{n-k}g(H)\ d\mu(H).
\end{equation}
When $k=1$ we get the original class of intersection bodies introduced by Lutwak in \cite{L}.

For a star body $K$ in $\R^n$ and $1\le k <n,$ denote by 
$${\rm {o.v.r.}}(K,{\cal{BP}}_k^n) = \inf \left\{ \left( \frac {|D|}{|K|}\right)^{1/n}:\ K\subset D,\ D\in {\cal{BP}}_k^n \right\}$$
the outer volume ratio distance from $K$ to the class ${\cal{BP}}_k^n.$ This quantity is directly related
to the isomorphic Busemann-Petty problem.

\begin{theorem}\label{ovr-bp} Suppose that $1\le k \le n-1,$ and  $K,L$ are
origin-symmetric star bodies in $\R^n$ such that
$$|K\cap H|\le |L\cap H|,\qquad \forall H\in Gr_{n-k}.$$
Then
$$|K|^{\frac{n-k}n}\le \left({\rm {o.v.r.}}(K,{\cal{BP}}_k^n)\right)^k|L|^{\frac{n-k}n}.$$
\end{theorem}

\pf Let $s>{\rm {o.v.r.}}(K,{\cal{BP}}_k^n),$ then there exists a star body
$D\in {\cal{BP}}_k^n$ such that $K\subset D$ and
\begin{equation}\label{ovr}
|D|^{\frac 1n}\le s |K|^{\frac 1n}.
\end{equation}
Let $\mu$ be the measure on $Gr_{n-k}$ corresponding to $D$ by definition (\ref{genint}).

By (\ref{polar-ld}), the condition $|K\cap H|\le |L\cap H|$ can be written as
$$R_{n-k}(\|\cdot\|_K^{-n+k})(H)\le R_{n-k}(\|\cdot\|_L^{-n+k})(H),\qquad \forall H\in Gr_{n-k}.$$
Integrating this inequality over $Gr_{n-k}$ with respect to the measure $\mu$ and using
(\ref{genint}) we get 
\begin{equation}\label{eq2}
\int_{S^{n-1}} \|x\|_D^{-k}\|x\|_K^{-n+k} dx \le \int_{S^{n-1}} \|x\|_D^{-k}\|x\|_L^{-n+k} dx.
\end{equation}
Since $K\subset D,$ we have $\|x\|_K^{-k}\le \|x\|_D^{-k},$ so the left-hand side of (\ref{eq2})
can be estimated from below by
$$\int_{S^{n-1}}\|x\|_K^{-n}dx = n|K|.$$
By H\"older's inequality and (\ref{ovr}), the right-hand side of (\ref{eq2}) can be estimated
from above by
$$\left(\int_{S^{n-1}} \|x\|_D^{-n} dx\right)^{\frac kn}\left(\int_{S^{n-1}} \|x\|_L^{-n} dx\right)^{\frac {n-k}n}=
n|D|^{\frac kn} |L|^{\frac {n-k}n}$$$$ \le n s^k
|K|^{\frac kn} |L|^{\frac{n-k}n}.$$
Combining these estimates and sending $s$ to ${\rm {o.v.r.}}(K,{\cal{BP}}_k^n),$ we get the result.
\endpf

The outer volume ratio distance from a general convex body to the class of
generalized $k$-intersection bodies was estimated in \cite{KPZ}.
\begin{pr} \label{kpz} (\cite[Theorem 1.1]{KPZ}) Let $K$ be an origin-symmetric convex 
body in $\R^n,$ and let $1\le k \le n-1.$ Then
$${\rm o.v.r.}(K,{\cal{BP}}_k^n) \le C \sqrt{\frac{n}{k}}\left(\log\left(\frac{en}{k}\right)\right)^{3/2},$$
where $C$ is an absolute constant.
\end{pr}
\smallbreak
\noindent {\bf Remark 1.} In \cite[Theorem 1.1]{KPZ}, the result was formulated with the logarithmic
term raised to the power 1/2 instead of 3/2, due to a mistake. The correction was made in \cite{K5}.
\medbreak
Our main result immediately follows.
\smallbreak
\noindent{\bf Proof of Theorem \ref{ildbp-th}.} Since $n/k < 1/\lambda,$ by Proposition \ref{kpz}, 
$${\rm o.v.r.}(K,{\cal{BP}}_k^n) < C \sqrt{\frac{(1-\log\lambda)^3}{\lambda}}.$$
The result follows from Theorem 2. \qed
\medbreak
\noindent{\bf Remark 2.} Theorem \ref{ildbp-th} implies Proposition \ref{ldsp} via a simple argument, similar
to the one employed in \cite{MP} for hyperplane sections. Indeed, assume that
$\lambda$ and $k$ are as in Theorem \ref{ildbp-th}, and suppose that
\begin{equation}\label{voleq}
|K|^{\frac{n-k}n} = C(\lambda)^k |B_2^n|^{\frac{n-k}n},
\end{equation}
where
$$C(\lambda)=C\sqrt{\frac{(1-\log \lambda)^3}{\lambda}},$$
$B_2^n$ is the unit Euclidean ball in $\R^n,$ and $K$ is an origin-symmetric convex body in $\R^n.$ 
Then, by Theorem \ref{ildbp-th},  it is not possible that $$|K\cap H|<|B_2^n\cap H|=|B_2^{n-k}|$$ for
all $H\in Gr_{n-k},$ so
$$\max_{H\in Gr_{n-k}} |K\cap H|\ge |B_2^{n-k}|.$$
Dividing both sides by equal numbers from (\ref{voleq}), we get
$$\frac{\max_{H\in Gr_{n-k}} |K\cap H|}{|K|^{\frac{n-k}n}} \ge 
\frac{|B_2^{n-k}|}{C(\lambda)^k |B_2^n|^{\frac{n-k}n}}.$$
By homogeneity, the condition (\ref{voleq}) can be dropped, and the latter inequality\
holds for arbitrary origin-symmetric convex $K.$ Now Proposition \ref{ldsp} follows from 
$$\frac{|B_2^n|^{\frac {n-k}n}}{|B_2^{n-k}|}\in (e^{-k/2},1);$$
see for example \cite[Lemma 2.1]{KL}.
\bigbreak
Finally, we mention the following fact which can be proved by applying Theorem \ref{ildbp-th} twice.
\begin{co} Suppose that $0<\lambda<1,$ $k>\lambda n,$ $0<c_1<c_2,$ and $K,L$ are origin-symmetric
convex bodies in $\R^n$ such that
$$c_1\le \frac{|L\cap H|}{|K\cap H|} \le c_2,\qquad \forall H\in Gr_{n-k}.$$
Then
$$\frac{c_1}{C^k\left(\sqrt{\frac{(1-\log \lambda)^3}{\lambda}}\right)^k}\le
\frac{|L|^{\frac{n-k}n}}{|K|^{\frac{n-k}n}} \le  c_2 C^k\left(\sqrt{\frac{(1-\log \lambda)^3}{\lambda}}\right)^k,$$
where $C$ is an absolute constant.
\end{co}

\bigbreak
{\bf Acknowledgement.} I wish to thank the US National Science Foundation for support through 
grant DMS-1265155.

\end{document}